\documentclass[11pt,oneside]{article}
\usepackage{psfrag}
\usepackage{epsfig}
\usepackage{graphicx}
\usepackage{amsmath}
\usepackage{amssymb}
\usepackage{amsthm}
\usepackage{color}
\usepackage{cite}
\usepackage{enumerate}
\usepackage[latin1]{inputenc}
\usepackage{multirow}
\usepackage{longtable}
\usepackage{multirow}

\newtheorem{theorem}{Theorem}

\begin{document}


\title{Semiovals in $\mathrm{PG}(2,8)$ and $\mathrm{PG}(2,9)$}
\date{\today}
\author{Daniele Bartoli, Stefano Marcugini, Fernanda Pambianco}

\maketitle



\begin{abstract}
The classification of all semiovals and blocking semiovals in $\mathrm{PG}(2,8)$ and in $\mathrm{PG}(2,9)$ of size less than $17$ is determined. Also, some information on the stabilizer groups and the intersection sizes with lines is given.
\end{abstract}

\emph{Keywords -} Semiovals, Blocking semiovals, Projective planes.


\section{Introduction}
Let $\Pi_q$ be a projective plane of order $q$. A {\it semioval} $\mathcal{S}$ in $\Pi_q$ is a non-empty pointset  with the property that for every point $P \in \mathcal{S}$ there exists a unique line $t_P$ such that $\mathcal{S}\cap t_P= \{P\}$ (see for example \cite{Kiss2007,CK2012}). This line is called the tangent line to $\mathcal{S}$ at $P$.

The classical examples of semiovals arise from polarities (ovals and unitals), and from the theory of blocking sets (the vertexless triangle). The semiovals are interesting objects in their own right, but the study of semiovals is also motivated by their applications to cryptography. Batten  constructed in \cite{Batten2000} an effective message sending scenario which uses determining sets. She showed that blocking semiovals are a particular type of determining sets in projective planes.

A {\it blocking semioval} (see \cite{Batten2000,Dover2000b,Dover2012}) is a semioval $\mathcal{S}$ such that every line of $\mathrm{PG}(2,q)$ contains at least one point of $\mathcal{S}$ and at least one point which is not in $\mathcal{S}$. A blocking semioval existing in every projective plane of order $q>2$ is the vertexless triangle.

In the last years the interest and research on the fundamental problem of determining the spectrum of the values for which there exists a given subconfiguration of points in $\mathrm{PG}(n,q)$ have increased considerably (see for example \cite{BDFMP2012,BMP2013b,BE1999,CS2012,DFMP2009,HS2001,KMP2010,MMP2005,MMP2007,PS2008}).

For $q\leq 9$, $q$ odd, the spectrum of semiovals was determined by Lisonek in \cite{Lisonek1994} by exhaustive computer search.

It is known that if $\mathcal{S}$ is a semioval in $\Pi _q$ then $q+1\leq |\mathcal{S} | \leq q\sqrt{q}+1$ and both bounds are sharp \cite{Hubaut1970,Thas1974}; the extremes occur when $\mathcal{S}$ is an oval or a unital.

The non-existence of semiovals of size $q+2$ in $\mathrm{PG}(2,q)$, $q\neq 7$ derives from  \cite[Theorem 4]{Blokhuis1991}.

A semioval $\mathcal{S}$ is called \emph{regular} if every non-tangent line of $\mathrm{PG}(2,q)$ intersects $\mathcal{S}$ in $0$ or in a constant number $a$ of points. G\'acs in \cite{Gacs2006} proved that if $\mathcal{S}$ is regular then $\mathcal{S}$ is a unital ($a=\sqrt{q}+1$) or an oval ($a=2$).

In this paper the full classification of semiovals and blocking semiovals in $\mathrm{PG}(2,8)$ and in $\mathrm{PG}(2,9)$ of size less than $17$ is given.

\section{Algorithm}\label{sect:A}
The algorithm used is a modification of the one presented in \cite{MMP2007}.

In this case, the algorithm works on \emph{admissible sets}, that is sets such that each point lies on at least one tangent line, instead of working on partial solutions. In fact, the property of being a semioval is not an  hereditary feature, that is a feature conserved by all its subsets, so the weaker hereditary feature of being an admissible set has been used. It is weaker in the sense that it allows to prune  few branches of the search space with respect to the cases when considering arcs and $3$-arcs. This and the fact that semiovals are in general longer than arcs and $3$-arcs make the problem computationally harder than the ones faced in \cite{MMP2005,MMP2007}.

Note also that, in general, not all the admissible sets can be extended to semiovals.

The exhaustive search has been feasible because projective properties among admissible sets have been exploited to avoid obtaining too many isomorphic copies of the same solution semioval and to avoid searching through parts of the search space isomorphic to previously searched ones.

The algorithm starts constructing a tree structure containing a representative of each class of non-equivalent admissible sets of size less than or equal to a fixed threshold $h$. If the threshold $h$ were equal to the actual size of the sought semiovals, the algorithm would be orderly, that is capable of constructing each goal configuration exactly once \cite{Royle1996}.

However, in the present case, the construction of the tree with the threshold $h$ equal to the size of the sought semiovals would have been too space and time consuming. For this reason a hybrid approach has been adopted. The  non-equivalent admissible sets of size $h$ obtained  have been extended using a backtracking algorithm, trying to determine semiovals of the desired size. In the backtracking phase, the information obtained during the classification of the admissible sets has been further exploited to prune the search tree. In fact, the points that would have given admissible sets equivalent to the ones already obtained have been excluded from the backtracking steps.

A simple parallelization technique, based on data distribution, has been used to divide the load of the computation in a multiprocessor computer. In the searches,  a $3.3$ Ghz Intel Exacore with 16 Gb of memory has been used.

\section{Semiovals in $\mathrm{PG}(2,8)$}\label{sect:PG2_8}
The spectrum of semiovals in $PG(2,8)$ has been determined in \cite[Theorem 2.1]{KMP2008} by computational methods.
\begin{theorem}
There exists a semioval of size $k$ in $PG(2,8)$ if and only if $k \in \{9,12-23\} $
\end{theorem}
Tables \ref{PG2_8}, \ref{PG2_8_int_sizes1}, \ref{PG2_8_int_sizes2}, and \ref{PG2_8gruppi} summarize the results obtained using the algorithm described in Section \ref{sect:A}. For every size, the numbers of non-equivalent semiovals and of blocking semiovals up to collineations in $P\Gamma L (3,8)$ are presented. For each size the spectrum of the possible intersection size with lines is also given: $\{a_1,\ldots,a_{k}\}$ indicates that there exist $i$-secants for each $i \in \{a_1,\ldots,a_{k}\}$ and the exponent indicates the number of non-equivalent examples having a particular spectrum. All the possible spectra of intersection size with lines are listed in detail in the subsequent tables, where $\ell_{i}$ indicates the number of $i$-secants.\\

\begin{center}
\begin{table}
\caption{Non-equivalent semiovals $\mathcal{S}$ in $PG(2,8)$}\label{PG2_8}
\begin{center}
\begin{tabular}{|c|c|c|c|}
\hline
{\bf Size}&Semiovals&Blocking&Intersection sizes with lines\\
&&semiovals&\\
\hline
\hline
{\bf 9}&$2$&$0$&$\{0,1,2\}^2$\\
\hline
{\bf 12}&$4$&$0$&$\{0,1,2,3\}^4$\\
\hline
{\bf 13}&$2$&$0$&$\{0,1,2,3,4\}^2$\\
\hline
{\bf 14}&$6$&$0$&$\{0,1,2,3,4\}^5, \{0,1,2,7\}^1$\\
\hline
{\bf 15}&$98$&$0$&\begin{tabular}{c} $\{0,1,2,3,4\}^{73}, \{0,1,2,3,4,5\}^{22}$ \\ $\{0, 1, 2, 3, 5, 6\}^{2}, \{ 0, 1, 2, 3, 5\}^{1}$\end{tabular}\\
\hline
{\bf 16}&$435$&$0$&\begin{tabular}{c} $\{0, 1, 2, 3, 4, 5, 6\}^{8}, \{0, 1, 2, 3, 4, 6\}^{7}$ \\ $\{0, 1, 2, 3, 4\}^{183}, \{  0, 1, 2, 3, 7\}^{1}$\\ $\{0, 1, 2, 3, 4, 5\}^{236}$ \\ \end{tabular}\\
\hline
{\bf 17}&$1064$&$0$&\begin{tabular}{c} $\{0, 1, 2, 3, 4, 5, 6\}^{23}, \{0, 1, 2, 3, 4, 6\}^{3}$ \\ $\{0, 1, 2, 3, 4\}^{226}, \{  0, 1, 2, 3, 7\}^{2}$\\ $\{0, 1, 2, 3, 4, 5\}^{810}$ \\ \end{tabular}\\
\hline
{\bf 18}&$1171$&$0$&\begin{tabular}{c} $\{0, 1, 2, 3, 4, 5, 6\}^{69}, \{0, 1, 2, 3, 4, 6\}^{5}$ \\ $\{0, 1, 2, 3, 4\}^{77}, \{0, 1, 2, 3, 6\}^{1}$\\ $\{0, 1, 2, 3, 4, 5\}^{1015}, \{0, 1, 2, 3, 5, 6 \}^{1}$ \\ $\{0, 1, 2, 3, 5\}^{1}, \{0, 1, 2, 3, 4, 7\}^{2}$\\ \end{tabular}\\
\hline
{\bf 19}&$884$&$2$&\begin{tabular}{c} $\{0, 1, 2, 3, 4, 5, 6\}^{111}, \{0, 1, 2, 3, 4, 6\}^{4}$\\  $\{0, 1, 2, 3, 4\}^{44}, \{0, 1, 2, 3, 4, 5\}^{722}$\\ $\{0, 1, 2, 3, 5, 7\}^{1}, \{1, 2, 3, 4, 5, 6\}^{2}$\\ \end{tabular}\\
\hline
{\bf 20}&$340$&$27$&\begin{tabular}{c}
    $\{ 0, 1, 2, 3, 4, 5, 6 \}^{29},\{ 0, 1, 2, 3, 6, 7 \}^{1}$\\
    $\{ 0, 1, 2, 3, 4 \}^{32},\{ 0, 1, 2, 3, 4, 5 \}^{251}$\\
    $\{ 1, 2, 3, 4, 5, 6 \}^{11},\{ 1, 2, 3, 4, 6, 7 \}^{3}$\\
    $\{ 1, 2, 3, 4, 5, 6, 7 \}^{1},\{ 1, 2, 3, 4, 5 \}^{12}$\\
\end{tabular}\\
\hline
{\bf 21}&$34$&$21$&
\begin{tabular}{c}
 $\{ 1, 3, 7 \}^{1}, \{ 0, 1, 2, 3, 4 \}^{8}$\\
    $\{ 0, 1, 2, 3, 4, 5 \}^{4},    \{ 1, 2, 3, 4, 6 \}^{1}$\\
    $\{ 1, 2, 3, 4, 5, 6 \}^{2},    \{ 0, 1, 3, 4 \}^{1}$\\
    $\{ 1, 2, 3, 4, 5 \}^{17}$\\
\end{tabular}\\
\hline
{\bf 22}&$1$&$0$&$\{0, 1, 2, 3, 4\}^1$\\
\hline
{\bf 23}&$1$&$1$&$\{1, 2, 3, 4\}^1$\\
\hline
\end{tabular}
\end{center}
\end{table}
\end{center}

\begin{center}
\begin{table}
\caption{Semiovals $\mathcal{S}$ in $PG(2,8)$: intersections with lines (Part 1)}\label{PG2_8_int_sizes1}
\begin{center}
{\footnotesize
\tabcolsep=1mm
\begin{tabular}{|c|c|c|c|c|c|c|c|c|c||c|c|c|c|c|c|c|c|c|c|}
\hline
{\bf Size}&$\ell_{0}$&$\ell_{1}$&$\ell_{2}$&$\ell_{3}$&$\ell_{4}$&$\ell_{5}$&$\ell_{6}$&$\ell_{7}$&Examples&{\bf Size}&$\ell_{0}$&$\ell_{1}$&$\ell_{2}$&$\ell_{3}$&$\ell_{4}$&$\ell_{5}$&$\ell_{6}$&$\ell_{7}$&Examples\\
\hline
\hline
\multirow{1}{*}{\bf 9}
&28& 9   & 36  &       &       &    &    &     & 2 &\multirow{9}{*}{\bf 17}&6 & 17 & 22 & 24 & 2 & & 2 & & 1 \\\cline{12-20}
\cline{1-10}
\cline{1-10}
\multirow{1}{*}{\bf 12}
&19& 12 & 30  & 12 &       &    &    &     & 4&&5 & 17 & 31 & 9 & 8 & 3 & & & 1 \\\cline{12-20}
\cline{1-10}
\cline{1-10}
\multirow{1}{*}{\bf 13}
&16& 13 & 30  & 12 &  2   &    &    &     & 2&&5 & 17 & 30 & 12 & 5 & 4 & & & 4 \\\cline{12-20}
\cline{1-10}
\cline{1-10}
\multirow{2}{*}{\bf 14}
&14& 14 & 25  & 18 &  2   &    &    &     & 5&&5 & 17 & 29 & 15 & 2 & 5 & & & 2 \\\cline{12-20}
\cline{2-10}
& 8 & 14 & 49  &       &       &    &    & 2  & 1&&5 & 17 & 29 & 14 & 5 & 2 & 1 & & 4 \\\cline{12-20}
\cline{1-10}
\cline{1-10}
\multirow{9}{*}{\bf 15}
&12 & 15 & 21 & 22 & 3 & & && 19 &&5 & 17 & 28 & 17 & 2 & 3 & 1 & & 1 \\\cline{12-20}
\cline{2-10}
&11 & 15 & 27 & 14 & 6 & & && 53 &&5 & 17 & 28 & 16 & 5 & & 2 & & 2 \\\cline{12-20}
\cline{2-10}
&11 & 15 & 26 & 17 & 3 & 1 & && 13 &&4 & 17 & 32 & 14 & 2 & 2 & 2 & & 1 \\\cline{12-20}
\cline{2-10}
&10 & 15 & 33 & 6 & 9 & & && 1 &&4 & 17 & 28 & 22 & & & & 2 & 2 \\ \cline{11-20}\cline{11-20}
\cline{2-10}
&10 & 15 & 32 & 9 & 6 & 1 & && 2 &\multirow{35}{*}{\bf 18}&7 & 18 & 9 & 30 & 9 & & & & 3\\ \cline{12-20}
\cline{2-10}
&10 & 15 & 31 & 12 & 3 & 2 & && 6&&7 & 18 & 9 & 30 & 9 & & & & 3  \\ \cline{12-20}
\cline{2-10}
&10 & 15 & 30 & 15 & & 3 & && 1 &&6 & 18 & 15 & 22 & 12 & & & & 60\\ \cline{12-20}
\cline{2-10}
&9 & 15 & 36 & 7 & 3 & 3 & && 1 &&6 & 18 & 14 & 25 & 9 & 1 & & & 56\\ \cline{12-20}
\cline{2-10}
&9 & 15 & 34 & 12 & & 2 & 1 && 2 &&6 & 18 & 13 & 28 & 6 & 2 & & & 15 \\\cline{12-20}
\cline{1-10}
\cline{1-10}
\multirow{19}{*}{\bf 16}
&11 & 16 & 12 & 32 & 2 & & & & 1 &&5 & 18 & 21 & 14 & 15 & & & & 14\\\cline{12-20}
\cline{2-10}
&10 & 16 & 18 & 24 & 5 & & & & 20 &&5 & 18 & 20 & 17 & 12 & 1 & & & 188\\\cline{12-20}
\cline{2-10}
&9 & 16 & 24 & 16 & 8 & & & & 157 &&5 & 18 & 19 & 20 & 9 & 2 & & & 368\\\cline{12-20}
\cline{2-10}
&9 & 16 & 23 & 19 & 5 & 1 & & & 113 &&5 & 18 & 18 & 23 & 6 & 3 & & & 70\\\cline{12-20}
\cline{2-10}
&9 & 16 & 22 & 22 & 2 & 2 & & & 2 &&5 & 18 & 18 & 22 & 9 & & 1 & & 6\\\cline{12-20}
\cline{2-10}
&8 & 16 & 30 & 8 & 11 & & & & 5 &&5 & 18 & 17 & 25 & 6 & 1 & 1 & & 7\\\cline{12-20}
\cline{2-10}
&8 & 16 & 29 & 11 & 8 & 1 & & & 48 &&4 & 18 & 26 & 9 & 15 & 1 & & & 1\\\cline{12-20}
\cline{2-10}
&8 & 16 & 28 & 14 & 5 & 2 & & & 63 &&4 & 18 & 25 & 12 & 12 & 2 & & & 29\\\cline{12-20}
\cline{2-10}
&8 & 16 & 27 & 17 & 2 & 3 & & & 8 &&4 & 18 & 24 & 15 & 9 & 3 & & & 161\\\cline{12-20}
\cline{2-10}
&8 & 16 & 27 & 16 & 5 & & 1 & & 2 &&4 & 18 & 24 & 14 & 12 & & 1 & & 3\\\cline{12-20}
\cline{2-10}
&8 & 16 & 26 & 19 & 2 & 1 & 1 & & 1 &&4 & 18 & 23 & 18 & 6 & 4 & & & 85\\\cline{12-20}
\cline{2-10}
&7 & 16 & 33 & 9 & 5 & 3 & & & 1 &&4 & 18 & 23 & 17 & 9 & 1 & 1 & & 13\\\cline{12-20}
\cline{2-10}
&7 & 16 & 33 & 8 & 8 & & 1 & & 2 &&4 & 18 & 22 & 21 & 3 & 5 & & & 17\\\cline{12-20}
\cline{2-10}
&7 & 16 & 32 & 12 & 2 & 4 & & & 1 &&4 & 18 & 22 & 20 & 6 & 2 & 1 & & 11\\\cline{12-20}
\cline{2-10}
&7 & 16 & 32 & 11 & 5 & 1 & 1 & & 5 &&4 & 18 & 21 & 23 & 3 & 3 & 1 & & 1\\\cline{12-20}
\cline{2-10}
&7 & 16 & 31 & 14 & 2 & 2 & 1 & & 1 &&4 & 18 & 21 & 22 & 6 & & 2 & & 3\\\cline{12-20}
\cline{2-10}
&7 & 16 & 30 & 16 & 2 & & 2 & & 3 &&4 & 18 & 18 & 30 & & & 3 & & 1\\\cline{12-20}
\cline{2-10}
&6 & 16 & 36 & 9 & 2 & 3 & 1 & & 1 &&3 & 18 & 30 & 7 & 12 & 3 & & & 1\\\cline{12-20}
\cline{2-10}
&5 & 16 & 36 & 14 & & & & 2 & 1 &&3 & 18 & 29 & 10 & 9 & 4 & & & 2\\\cline{12-20}
\cline{1-10}
\cline{1-10}
\multirow{13}{*}{\bf 17}
&8 & 17 & 16 & 24 & 8 & & & & 59 &&3 & 18 & 28 & 13 & 6 & 5 & & & 12\\\cline{12-20}
\cline{2-10}
&8 & 17 & 15 & 27 & 5 & 1 & & & 10 &&3 & 18 & 28 & 12 & 9 & 2 & 1 & & 11\\\cline{12-20}
\cline{2-10}
&7 & 17 & 22 & 16 & 11 & & & & 165 &&3 & 18 & 27 & 16 & 3 & 6 & & & 1\\\cline{12-20}
\cline{2-10}
&7 & 17 & 21 & 19 & 8 & 1 & & & 388 &&3 & 18 & 27 & 15 & 6 & 3 & 1 & & 19\\\cline{12-20}
\cline{2-10}
&7 & 17 & 20 & 22 & 5 & 2 & & & 164 &&3 & 18 & 27 & 14 & 9 & & 2 & & 2\\\cline{12-20}
\cline{2-10}
&7 & 17 & 19 & 25 & 2 & 3 & & & 2 &&3 & 18 & 26 & 18 & 3 & 4 & 1 & & 3\\\cline{12-20}
\cline{2-10}
&6 & 17 & 28 & 8 & 14 & & & & 2 &&3 & 18 & 26 & 17 & 6 & 1 & 2 & & 2\\\cline{12-20}
\cline{2-10}
&6 & 17 & 27 & 11 & 11 & 1 & & & 22&&3 & 18 & 25 & 21 & & 5 & 1 & & 1 \\\cline{12-20}
\cline{2-10}
&6 & 17 & 26 & 14 & 8 & 2 & & & 139 &&3 & 18 & 25 & 20 & 3 & 2 & 2 & & 1\\\cline{12-20}
\cline{2-10}
&6 & 17 & 25 & 17 & 5 & 3 & & & 71 &&3 & 18 & 21 & 28 & 1 & & & 2 & 2 \\\cline{12-20}
\cline{2-10}
&6 & 17 & 25 & 16 & 8 & & 1 & & 7 &&2 & 18 & 31 & 12 & 6 & 2 & 2 & & 1\\\cline{12-20}
\cline{2-10}
&6 & 17 & 24 & 19 & 5 & 1 & 1 & & 15 &&1 & 18 & 36 & 9 & & 9 & & & 1 \\\cline{11-20}\cline{11-20}

\cline{2-10}
&6 & 17 & 23 & 22 & 2 & 2 & 1 & & 2 &\multirow{1}{*}{\bf 19}&5 & 19 & 9 & 26 & 14 & & & & 11\\\cline{11-20}
\cline{1-10}
\end{tabular}
}
\end{center}
\end{table}
\end{center}

\begin{center}
\begin{table}
\caption{Semiovals $\mathcal{S}$ in $PG(2,8)$: intersections with lines (Part 2)}\label{PG2_8_int_sizes2}
\begin{center}
{\footnotesize
\tabcolsep=1mm
\begin{tabular}{|c|c|c|c|c|c|c|c|c|c||c|c|c|c|c|c|c|c|c|c|}
\hline
Size&$\ell_{0}$&$\ell_{1}$&$\ell_{2}$&$\ell_{3}$&$\ell_{4}$&$\ell_{5}$&$\ell_{6}$&$\ell_{7}$&Examples&Size&$\ell_{0}$&$\ell_{1}$&$\ell_{2}$&$\ell_{3}$&$\ell_{4}$&$\ell_{5}$&$\ell_{6}$&$\ell_{7}$&Examples\\
\hline
\hline
\multirow{36}{*}{\bf 19}
&5 & 19 & 8 & 29 & 11 & 1 & & & 1&\multirow{28}{*}{\bf 20}&2 & 20 & 12 & 24 & 11 & 4 & & & 24 \\\cline{12-20}
\cline{2-10}
&4 & 19 & 15 & 18 & 17 & & & & 32 &&2 & 20 & 12 & 23 & 14 & 1 & 1 & & 1\\\cline{12-20}
\cline{2-10}
&4 & 19 & 14 & 21 & 14 & 1 & & & 50&&2 & 20 & 11 & 27 & 8 & 5 & & & 1 \\\cline{12-20}
\cline{2-10}
&4 & 19 & 13 & 24 & 11 & 2 & & & 49 &&1 & 20 & 21 & 7 & 23 & 1 & & & 1 \\\cline{12-20}
\cline{2-10}
&4 & 19 & 12 & 27 & 8 & 3 & & & 16 &&1 & 20 & 20 & 10 & 20 & 2 & & & 3\\\cline{12-20}
\cline{2-10}
&3 & 19 & 21 & 10 & 20 & & & & 1 &&1 & 20 & 19 & 13 & 17 & 3 & & & 19\\\cline{12-20}
\cline{2-10}
&3 & 19 & 20 & 13 & 17 & 1 & & & 5&&1 & 20 & 18 & 16 & 14 & 4 & & & 57 \\\cline{12-20}
\cline{2-10}
&3 & 19 & 19 & 16 & 14 & 2 & & & 113&&1 & 20 & 18 & 15 & 17 & 1 & 1 & & 3 \\\cline{12-20}
\cline{2-10}
&3 & 19 & 18 & 19 & 11 & 3 & & & 219 &&1 & 20 & 17 & 19 & 11 & 5 & & & 51\\\cline{12-20}
\cline{2-10}
&3 & 19 & 18 & 18 & 14 & & 1 & & 2 &&1 & 20 & 17 & 18 & 14 & 2 & 1 & & 10\\\cline{12-20}
\cline{2-10}
&3 & 19 & 17 & 22 & 8 & 4 & & & 100 &&1 & 20 & 16 & 22 & 8 & 6 & & & 13\\\cline{12-20}
\cline{2-10}
&3 & 19 & 17 & 21 & 11 & 1 & 1 & & 6 &&1 & 20 & 16 & 21 & 11 & 3 & 1 & & 9\\\cline{12-20}
\cline{2-10}
&3 & 19 & 16 & 25 & 5 & 5 & & & 10 &&1 & 20 & 15 & 25 & 5 & 7 & & & 4 \\\cline{12-20}
\cline{2-10}
&3 & 19 & 16 & 24 & 8 & 2 & 1 & & 5 &&1 & 20 & 15 & 24 & 8 & 4 & 1 & & 5 \\\cline{12-20}
\cline{2-10}
&3 & 19 & 15 & 27 & 5 & 3 & 1 & & 1 &&1 & 20 & 14 & 27 & 5 & 5 & 1 & & 1\\\cline{12-20}
\cline{2-10}
&2 & 19 & 25 & 8 & 17 & 2 & & & 1 &&1 & 20 & 7 & 42 & & & 1 & 2 & 1\\\cline{12-20}
\cline{2-10}
&2 & 19 & 24 & 11 & 14 & 3 & & & 8 && & 20 & 24 & 7 & 20 & 1 & 1 & & 1\\\cline{12-20}
\cline{2-10}
&2 & 19 & 23 & 14 & 11 & 4 & & & 63&& & 20 & 23 & 11 & 14 & 5 & & & 2 \\\cline{12-20}
\cline{2-10}
&2 & 19 & 23 & 13 & 14 & 1 & 1 & & 5 && & 20 & 23 & 10 & 17 & 2 & 1 & & 1\\\cline{12-20}
\cline{2-10}
&2 & 19 & 22 & 17 & 8 & 5 & & & 73 && & 20 & 22 & 14 & 11 & 6 & & & 4\\\cline{12-20}
\cline{2-10}
&2 & 19 & 22 & 16 & 11 & 2 & 1 & & 31&& & 20 & 22 & 13 & 14 & 3 & 1 & & 3 \\\cline{12-20}
\cline{2-10}
&2 & 19 & 21 & 20 & 5 & 6 & & & 11 && & 20 & 21 & 17 & 8 & 7 & & & 6\\\cline{12-20}
\cline{2-10}
&2 & 19 & 21 & 19 & 8 & 3 & 1 & & 34& & & 20 & 21 & 16 & 11 & 4 & 1 & & 2 \\\cline{12-20}
\cline{2-10}
&2 & 19 & 21 & 18 & 11 & & 2 & & 2 && & 20 & 21 & 15 & 14 & 1 & 2 & & 1\\\cline{12-20}
\cline{2-10}
&2 & 19 & 20 & 22 & 5 & 4 & 1 & & 8 && & 20 & 20 & 19 & 8 & 5 & 1 & & 1\\\cline{12-20}
\cline{2-10}
&2 & 19 & 20 & 21 & 8 & 1 & 2 & & 4&& & 20 & 20 & 18 & 11 & 2 & 2 & & 1 \\\cline{12-20}
\cline{2-10}
&2 & 19 & 14 & 35 & & 1 & & 2 & 1& & & 20 & 19 & 21 & 8 & 3 & 2 & & 1\\\cline{12-20}
\cline{2-10}
&1 & 19 & 28 & 9 & 11 & 5 & & & 1 && & 20 & 19 & 19 & 13 & & 1 & 1 & 3\\\cline{12-20}
\cline{2-10}
&1 & 19 & 27 & 12 & 8 & 6 & & & 1 && & 20 & 18 & 22 & 10 & 1 & 1 & 1 & 1 \\\cline{11-20} \cline{11-20}
\cline{2-10}
&1 & 19 & 27 & 11 & 11 & 3 & 1 & & 5 &\multirow{11}{*}{21}&3 & 21 & & 28 & 21 & & & & 1\\\cline{12-20}
\cline{2-10}
&1 & 19 & 26 & 15 & 5 & 7 & & & 1 &&2 & 21 & 6 & 20 & 24 & & & & 8\\\cline{12-20}
\cline{2-10}
&1 & 19 & 26 & 14 & 8 & 4 & 1 & & 4 &&1 & 21 & 9 & 21 & 18 & 3 & & & 4\\\cline{12-20}
\cline{2-10}
&1 & 19 & 25 & 17 & 5 & 5 & 1 & & 2 && & 21 & 15 & 13 & 21 & 3 & & & 1\\\cline{12-20}
\cline{2-10}
&1 & 19 & 25 & 16 & 8 & 2 & 2 & & 3 && & 21 & 15 & 12 & 24 & & 1 & & 1 \\\cline{12-20}
\cline{2-10}
&1 & 19 & 24 & 19 & 5 & 3 & 2 & & 3 && & 21 & 14 & 16 & 18 & 4 & & & 8\\\cline{12-20}
\cline{2-10}
& & 19 & 30 & 12 & 5 & 6 & 1 & & 2 && & 21 & 13 & 19 & 15 & 5 & & & 3\\\cline{12-20}
\cline{1-10}
\cline{1-10}
\multirow{9}{*}{\bf 20}&4 & 20 & 4 & 28 & 17 & & & & 3 &&& 21 & 12 & 22 & 12 & 6 & & & 5\\\cline{12-20}
\cline{2-10}
&3 & 20 & 10 & 20 & 20 & & & & 28 &&& 21 & 12 & 21 & 15 & 3 & 1 & & 1\\\cline{12-20}
\cline{2-10}
&3 & 20 & 8 & 26 & 14 & 2 & & & 2 && & 21 & 9 & 30 & 6 & 6 & 1 & & 1\\\cline{12-20}
\cline{2-10}
&3 & 20 & 7 & 29 & 11 & 3 & & & 1&& & 21 & & 49 & & & & 3 & 1 \\\cline{11-20}\cline{11-20}
\cline{2-10}
&2 & 20 & 16 & 12 & 23 & & & & 1 &\multirow{1}{*}{\bf 22}&1&22&3&18&29&&&&1\\\cline{11-20}\cline{11-20}
\cline{2-10}
&2 & 20 & 15 & 15 & 20 & 1 & & & 4 &\multirow{1}{*}{\bf 23}&&23&1&14&35&&&&1\\\cline{11-20}\cline{11-20}
\cline{2-10}
&2 & 20 & 14 & 18 & 17 & 2 & & & 31 \\
\cline{2-10}
&2 & 20 & 13 & 21 & 14 & 3 & & & 40 \\
\cline{1-10}
\end{tabular}
}
\end{center}
\end{table}
\end{center}

\begin{center}
\begin{table}
\caption{Semiovals $\mathcal{S}$ in $PG(2,8)$: stabilizer groups}\label{PG2_8gruppi}
\begin{center}
\begin{tabular}{|c||ccccc|}
\hline
{\bf Size}&&&&&\\
\hline
\hline
{\bf 9}& $G_{168}$: $1$&$G_{1512}$: $1$&&&\\
\hline
{\bf 12}& $\mathbb{Z}_3$: $2$&$\mathbb{Z}_{12}$: $1$&$G_{36}$ : $1$&&\\
\hline
{\bf 13}& $\mathbb{Z}_2$: $1$&$\mathbb{Z}_6$: $1$&&&\\
\hline
{\bf 14}& $\mathbb{Z}_1$: $2$&$\mathbb{Z}_3$: $2$&$\mathbb{Z}_6$: $1$&$G_{294}$: $1$&\\
\hline
{\bf 15}& $\mathbb{Z}_1$: $81$&$\mathbb{Z}_2$: $8$&$\mathbb{Z}_3$: $5$&$\mathbb{Z}_6$: $3$&$\mathcal{S}_3$: $1$\\
\hline
\multirow{2}*{{\bf 16}}& $\mathbb{Z}_1$: $412$&$\mathbb{Z}_2$: $7$&$\mathbb{Z}_3$: $12$&$\mathbb{Z}_4$: $1$&$\mathcal{D}_4$: $1$\\
&$\mathcal{Q}_6$: $1$&$\mathbb{Z}_{14}$: $1$&&&\\
\hline
\multirow{2}*{{\bf 17}}& $\mathbb{Z}_1$: $1014$&$\mathbb{Z}_2$: $25$&$\mathbb{Z}_3$: $21$&$\mathbb{Z}_2 \times \mathbb{Z}_2$: $2$&$\mathbb{Z}_{14}$: $1$\\
&$G_{21}$: $1$&&&&\\
\hline
\multirow{2}*{{\bf 18}}& $\mathbb{Z}_1$: $1133 $&$\mathbb{Z}_2$: $5$&$\mathbb{Z}_3$: $26$&$\mathbb{Z}_6$: $2$&$\mathbb{Z}_{3}\times \mathbb{Z}_3$: $1$\\
&$\mathbb{Z}_9$: $1$&$\mathbb{Z}_{14}$: $1$&$G_{18}$: $1$&$G_{21}$: $1$&\\
\hline
\multirow{1}*{{\bf 19}}& $\mathbb{Z}_1$: $851 $&$\mathbb{Z}_2$: $26$&$\mathbb{Z}_3$: $3$&$\mathbb{Z}_6$: $3$&$\mathbb{Z}_{14}$: $1$\\
\hline
\multirow{1}*{{\bf 20}}& $\mathbb{Z}_1$: $321 $&$\mathbb{Z}_3$: $18$&$G_{42}$: $1$&&\\
\hline
\multirow{1}*{{\bf 21}}& $\mathbb{Z}_1$: $26$&$\mathbb{Z}_3$: $6$&$G_{63}$: $1$&$G_{882}$: $1$&\\
\hline
\multirow{1}*{{\bf 22}}& $\mathbb{Z}_1$: $1$&&&&\\
\hline
\multirow{1}*{{\bf 23}}& $G_{21}$: $1$&&&&\\
\hline
\end{tabular}
\end{center}
\end{table}
\end{center}

\begin{center}
\begin{table}
\caption{Semiovals with a $7$-secant in $PG(2,8)$}\label{7_secant}
\begin{center}
\tabcolsep=0.75 mm
\begin{tabular}{|c|c|c|c|c|c|c|c|c|c|}
\hline
$\mathcal{S}$&$\ell_{0}$&$\ell_{1}$&$\ell_{2}$&$\ell_{3}$&$\ell_{4}$&$\ell_{5}$&$\ell_{6}$&$\ell_{7}$&$G$\\
\hline
\hline
\tabcolsep=0.75 mm
\begin{tabular}{cccccccccccccccccccc}
 1& 0& 0& 1& 0& 0& 1& 1& 1& 1& 1& 1& 1& 1& 1& 1& 1& 1& 1& 1\\
 0& 1& 0& 1& 1& 1& 0& 0& 1& 1& 3& 4& 5& 5& 5& 5& 5& 5& 6& 7\\
 0& 0& 1& 1& 2& 3& 3& 4& 2& 3& 3& 3& 0& 1& 2& 4& 5& 7& 7& 3\\
\end{tabular}&
0 &
20 &
19 &
19 &
13 &
0 &
1 &
1 &
$\mathbb{Z}_{1}$\\
\hline
\tabcolsep=0.75 mm
\begin{tabular}{cccccccccccccccccccc}
 1& 0& 0& 1& 0& 0& 1& 1& 1& 1& 1& 1& 1& 1& 1& 1& 1& 1& 1& 1\\
 0& 1& 0& 1& 1& 1& 0& 1& 1& 3& 4& 5& 5& 5& 5& 5& 5& 6& 7& 7\\
 0& 0& 1& 1& 2& 3& 3& 2& 3& 3& 3& 0& 1& 2& 4& 5& 7& 7& 1& 3\\
\end{tabular}&
0 &
20 &
18 &
22 &
10 &
1 &
1 &
1 &
$\mathbb{Z}_{1}$\\
\hline
\tabcolsep=0.75 mm
\begin{tabular}{cccccccccccccccccccc}
 1& 0& 0& 1& 0& 0& 1& 1& 1& 1& 1& 1& 1& 1& 1& 1& 1& 1& 1& 1\\
 0& 1& 0& 1& 1& 1& 0& 0& 1& 3& 4& 4& 5& 5& 5& 5& 5& 5& 6& 7\\
 0& 0& 1& 1& 2& 3& 3& 4& 3& 3& 3& 7& 0& 1& 2& 4& 5& 7& 7& 3\\
\end{tabular}&
0 &
20 &
19 &
19 &
13 &
0 &
1 &
1 &
$\mathbb{Z}_{3}$\\
\hline
\tabcolsep=0.75 mm
\begin{tabular}{cccccccccccccccccccc}
 1& 0& 0& 1& 0& 1& 1& 1& 1& 1& 1& 1& 1& 1& 1& 1& 1& 1& 1& 1\\
 0& 1& 0& 1& 1& 0& 1& 2& 3& 3& 5& 5& 5& 5& 5& 5& 6& 6& 7& 7\\
 0& 0& 1& 1& 2& 3& 2& 4& 3& 4& 0& 1& 2& 4& 5& 7& 2& 7& 1& 3\\
\end{tabular}&
0 &
20 &
19 &
19 &
13 &
0 &
1 &
1 &
$\mathbb{Z}_{3}$\\
\hline
\end{tabular}
\end{center}
\end{table}
\end{center}

\subsection{Semiovals in $\mathrm{PG}(2,9)$ of size less than $17$}
The spectrum of semiovals in $PG(2,9)$ has been determined in \cite{Lisonek1994}.
\begin{theorem}
There exists a semioval of size $k$ in $PG(2,9)$ if and only if $k \in \{ 10,12-28\} .$
\end{theorem}
We classified by computer search the semiovals $\mathcal{S}$ in $PG(2,9)$ of size less than or equal to $16$.
Tables \ref{PG2_9}, \ref{PG2_9_int_sizes}, and \ref{PG2_9gruppi} summarize the results obtained.

\begin{center}
\begin{table}
\caption{Non-equivalent semiovals $\mathcal{S}$ in $PG(2,9)$, with $|\mathcal{S}|\leq 16$}\label{PG2_9}
\begin{center}
\begin{tabular}{|c||c|c|c|c|c|c|}
\hline
Size&$10$&$12$&$13$&$14$&$15$&$16$\\
\hline
Semiovals&$1$&$1$&$1$&$3$&$26$&$113$\\

\hline
\end{tabular}
\end{center}
\end{table}
\end{center}

\begin{center}
\begin{table}
\caption{Semiovals $\mathcal{S}$ in $PG(2,9)$ with $|\mathcal{S}|\leq16$: intersections with lines}\label{PG2_9_int_sizes}
\begin{center}
\begin{tabular}{|c|c|c|c|c|c|c|c|c|c|c|}
\hline
{\bf Size}&$\ell_{0}$&$\ell_{1}$&$\ell_{2}$&$\ell_{3}$&$\ell_{4}$&$\ell_{5}$&$\ell_{6}$&$\ell_{7}$&$\ell_8$&Examples\\
\hline
\hline
\multirow{1}{*}{\bf 10}
&36 & 10 & 45 & &&&&&&1 \\
\hline
\hline

\multirow{1}{*}{\bf 12}
&28 & 12 & 48 & & 3 &&&&& 1 \\
\hline
\hline
\multirow{1}{*}{\bf 13}
&25 & 13 & 44 & 8 & & 1 &&&& 1 \\
\hline
\hline
\multirow{1}{*}{\bf 14}
&23 & 14 & 37 & 16 & 1 &&&&& 3 \\
\hline
\hline
\multirow{3}{*}{\bf 15}
&21 & 15 & 30 & 25 & & &&&&4 \\
\cline{2-11}
&20 & 15 & 36 & 17 & 3 & &&&&18 \\
\cline{2-11}
&19 & 15 & 42 & 9 & 6 & &&&&4 \\
\hline
\hline
\multirow{10}{*}{\bf 16}
&19 & 16 & 24 & 32 & & & & & & 2 \\\cline{2-11}
&18 & 16 & 30 & 24 & 3 & & & & & 39 \\\cline{2-11}
&17 & 16 & 36 & 16 & 6 & & & & & 51 \\\cline{2-11}
&17 & 16 & 35 & 19 & 3 & 1 & & & & 8 \\\cline{2-11}
&16 & 16 & 42 & 8 & 9 & & & & & 3 \\\cline{2-11}
&16 & 16 & 41 & 11 & 6 & 1 & & & & 2 \\\cline{2-11}
&16 & 16 & 40 & 14 & 3 & 2 & & & & 4 \\\cline{2-11}
&15 & 16 & 48 & & 12 & & & & & 2 \\\cline{2-11}
&15 & 16 & 45 & 8 & 6 & & 1 & & & 1 \\\cline{2-11}
&9 & 16 & 64 & & & & & & 2 & 1 \\
\hline
\end{tabular}
\end{center}
\end{table}
\end{center}

\begin{center}
\begin{table}
\caption{Semiovals $\mathcal{S}$ in $PG(2,9)$ with $|\mathcal{S}|\leq16$: stabilizer groups}\label{PG2_9gruppi}
\begin{center}
\begin{tabular}{|c||ccccc|}
\hline
{\bf Size}&&&&&\\
\hline
\hline
{\bf 10}& $G_{1440}$: $1$&&&&\\
\hline
{\bf 12}&$G_{192}$ : $1$&&&&\\
\hline
{\bf 13}&$G_{16}$ : $1$&&&&\\
\hline
{\bf 14}& $\mathbb{Z}_2$: $2$&$\mathbb{Z}_2\times \mathbb{Z}_2$: $1$&&&\\
\hline
\multirow{2}*{\bf 15}& $\mathbb{Z}_1$: $6$&$\mathbb{Z}_2$: $11$&$\mathbb{Z}_3$: $2$&$\mathbb{Z}_6$: $2$&$\mathcal{S}_3$: $2$\\
&$\mathcal{D}_4$: $1$&$\mathcal{Q}_6$: $1$&&&\\
\hline
\multirow{3}*{{\bf 16}}& $\mathbb{Z}_1$: $81$&$\mathbb{Z}_2$: $22$&$\mathbb{Z}_2 \times \mathbb{Z}_2$: $2$&$\mathbb{Z}_4$: $1$&$\mathbb{Z}_6$: $1$\\
&$\mathbb{Z}_2 \times \mathbb{Z}_4$: $1$&$\mathcal{D}_{4}$: $1$&$G_{32}$: $1$&$G_{48}$: $1$&$G_{192}$: $1$\\
&$G_{256}$: $1$&&&&\\
\hline
\end{tabular}
\end{center}
\end{table}
\end{center}

\end{document}